\documentclass[12pt]{article}%
\usepackage{amsmath,amssymb,fullpage}
\usepackage[applemac]{inputenc}
\usepackage{amsmath}
\usepackage{amsfonts}
\usepackage{amssymb}
\usepackage{graphicx}%
\providecommand{\U}[1]{\protect\rule{.1in}{.1in}}
\newtheorem{example}{Example}[section]
\newtheorem{theorem}[example]{Theorem}
\newtheorem{definition}[example]{Definition}
\newtheorem{lemma}[example]{Lemma}

\newtheorem{remark}[example]{Remark}

\numberwithin{equation}{section}

\def\1B{\text{1\!\!I}}

\begin{document}

\title{Mean-Field Delayed BSDEs with Jumps}
\author{Nacira Agram}
\date{10 January 2018}
\maketitle

\footnotetext[1]{Department of Mathematics, University of Oslo, P.O. Box 1053
Blindern, N--0316 Oslo, Norway and University of Biskra, Algeria. Email:
naciraa@math.uio.no.
\par
{}}

\footnotetext[2]{This research was carried out with support of the Norwegian
Research Council, within the research project Challenges in Stochastic
Control, Information and Applications (STOCONINF), project number 250768/F20.}

\paragraph{Abstract}

\noindent We establish sufficient conditions for the existence and uniqueness
of mean-field backward stochastic differential equations with time delayed
generator in the sense that at $t$, the generator may depend on previous
values up to a delay constant $\delta$ not on the hole past as in Delong and
Imkeller \cite{di},\cite{DI2}. For sufficiently small delay constant $\delta$
and for any finite time horizon, we get a unique solution.

\paragraph{MSC(2010):}

\noindent34F05, 60H30, 60G51, 60G57.

\paragraph{Keywords:}

\noindent Backward stochastic differential equations; Time delayed generator;
Mean-field; Lévy process; Poisson random measure.

\section{Introduction}

\noindent Backward stochastic differential equations (BSDEs) appear in their
linear form as an adjoint equation when Bismut \cite{b} was dealing with
stochastic optimal control. After that, this theory has been developed also
for the nonlinear case, we refer to namely to the seminal work by Pardoux and
Peng \cite{PP} and also to Pardoux \cite{P}, El Karoui \textit{et al}
\cite{EPQ}. The first work applying in finance was made by El Karoui
\textit{et al} \cite{EPQ}.

\noindent Given a driven Brownian motion $B$, a generator $f:\Omega
\times\lbrack0,T]\times%
%TCIMACRO{\U{211d} }%
%BeginExpansion
\mathbb{R}
%EndExpansion
^{2}\rightarrow%
%TCIMACRO{\U{211d} }%
%BeginExpansion
\mathbb{R}
%EndExpansion
$ and a terminal condition $\xi.$ Solving a BSDE consists in finding a couple
of processes $\left(  Y(t),Z(t)\right)  _{t\geq0}$ adapted to the considered
filtration (the Brownian one), such that, at time $t$, $\left(
Y(t),Z(t)\right)  _{t\geq0}$ satisfies the equation%
\begin{equation}
Y(t)=\xi+%
%TCIMACRO{\tint _{t}^{T}}%
%BeginExpansion
{\textstyle\int_{t}^{T}}
%EndExpansion
f(s,Y(s),Z(s))ds-%
%TCIMACRO{\tint _{t}^{T}}%
%BeginExpansion
{\textstyle\int_{t}^{T}}
%EndExpansion
Z(s)dB(s),\text{ }0\leq t\leq T. \label{n}%
\end{equation}

\noindent The crucial question now is which conditions should be satisfied by
the generator $f$ and the terminal value $\xi$ in order to get the existence
and the uniqueness of such a solution to $\left(  \ref{n}\right)  $.

\noindent In this paper we are interested in a generalisation of this above
BSDE, where at time $s$ the coefficient $f$ depends on past information and
the law of the solution process and we consider also the discountinuous case.
More precisely, we are interested in the Mean-Field Delayed BSDE (MF-DBSDE)
with jumps of the form%

\begin{equation}
\left\{
\begin{array}
[c]{ll}%
Y(t) & =\xi+%
%TCIMACRO{\tint _{t}^{T}}%
%BeginExpansion
{\textstyle\int_{t}^{T}}
%EndExpansion
f(s,Y_{s},Z_{s},K_{s}(\cdot),P_{(Y_{s},Z_{s},K_{s}(\cdot))})ds-%
%TCIMACRO{\tint _{t}^{T}}%
%BeginExpansion
{\textstyle\int_{t}^{T}}
%EndExpansion
Z(s)dB(s)\\
& \text{ \ \ \ \ \ \ \ \ \ \ \ \ }-%
%TCIMACRO{\tint _{t}^{T}}%
%BeginExpansion
{\textstyle\int_{t}^{T}}
%EndExpansion%
%TCIMACRO{\tint _{\mathbb{R}_{0}}}%
%BeginExpansion
{\textstyle\int_{\mathbb{R}_{0}}}
%EndExpansion
K(s,\zeta)\tilde{N}(ds,d\zeta),t\in\left[  0,T\right]  ,\\
Y(t) & =Y(0),Z(t)=0,K(t,\cdot)=0,t<0.
\end{array}
\right.  \label{1st}%
\end{equation}
Here, for a given delay constant $\delta>0$,
\[%
\begin{array}
[c]{ll}%
Y_{s} & =\left(  Y(s+r)\right)  _{r\in\left[  -\delta,0\right]  },\\
Z_{s} & =\left(  Z(s+r)\right)  _{r\in\left[  -\delta,0\right]  },\\
K_{s}(\cdot) & =\left(  K(s+r,\cdot)\right)  _{r\in\left[  -\delta,0\right]
},
\end{array}
\]
where information on the past of the solution process $(Y,Z,K)$ is considered.
We remark that in this case, i.e., when the dependence of the coefficient $f$
on past information is studied, the process $f(s,Y_{s},Z_{s},K_{s}%
(\cdot),P_{(Y_{s},Z_{s},K_{s}(\cdot))})$ is already $\mathbb{F}${-adapted, and
coincides }${ds}${ }${dP}${-a.e. with its optional projection. }

\noindent The case of DBSDEs in both continuous and discontinuous case have
been studied by Delong and Imkeller \cite{di} and \cite{DI2}. The authors
consider the coefficients depending on the full past of the path of the
solution. Following the same approach of Delong and Imkeller \cite{di}, Agram
and Røse \cite{AR2} obtained existence and uniqueness of MF-DBSDE where the
mean-field term considered is the expectation of the state and also we refer
to the paper by Ma and Liu \cite{ml}, in this framework.

\noindent The case when the mean-feld is represented by the law of the state
has been studied by Carmona and Delarue \cite{CD} and with discrete delay and
implicite terminal condition has been studied by Agram \cite{A}. The last two
papers, the probability measures are defined in the Wasserstein metric space
$\mathcal{P}_{2}$ and the Wasserstein distance $W_{2}$ is used.

\noindent In this note we study MF-BSDEs of type (\ref{1st}) driven by a
Brownian motion and a jumps and the filtration is that generated by both the
Brownian motion and the independent Poisson random measure\textbf{.} We
consider the Hilbert space of measures $\mathcal{M}$ introduced by Agram and
Øksendal \cite{AO2}, \cite{AO3}. The delay of the state processes is taking up
to a delay constant $\delta$ not of the hole past as in Delong and Imkeller
\cite{di} and this will help to get existence and uniqueness of the MF-DBSDE
(\ref{1st}) for any finite time horizon $T$ and any Lipschitz constant $C$ but
for sufficiently small delay constant $\delta$.

\noindent BSDEs with jumps have been studied by many authors, for exampe, we
refer to Tang and Li \cite{TL}, Barles \textit{et al} \cite{BBP}, Royer
\cite{R} Sulem and Quenez \cite{SQ} and Øksendal and Sulem \cite{OS}%
.\smallskip

\noindent These type of MF-DBSDE generalises the classical BSDE and have
turned out to be useful in various applications, namely in finance and in
stochastic control. While in finance the delay imposes in the modelling by the
fact that agents have often only time-delayed information, for example: If one
wants to find an investment strategy and an investment portfolio which
replicate a liability or meet a target which depends on the applied strategy
or the portfolio depends on its past values, then the DBSDEs are the best tool
to solve this financial problem. BSDEs with delay can also arise in portfolio
management problems, variable annuities, unit-linked products and
participating contract. For more details about applications of such equations,
we refer to Delong \cite{d}.

\noindent Next sections are devoted to the study of the MF-DBSDEs as follows:

\begin{itemize}
\item in section $2$, we introduce the adequates spaces of processes and we
fixe suitable assumptions on the driver and the terminal value.

\item in section $3$, we give a theorem on the existence and the uniqueness of
the solutions.
\end{itemize}

\section{Background}

\noindent Let $B(t),t\geq0$ be a $1$-dimentional Brownian motion, and
$\tilde{N}(dt,d\zeta),t\geq0$ be an independent compensated Poisson random
measure, with compensator $\nu(d\zeta)dt$, on a probability space
$(\Omega,\mathcal{F},P)$. Let $\mathbb{F=}(\mathcal{F}_{t})_{t\geq0}$ denote
the natural filtration satisfying the usual conditions of right continuity and
completeness, associated with $B$ and $N$. Let $\delta>0$, and extend the
filtration by letting $\mathcal{F}_{t}=\mathcal{F}_{0}$ for $t\in
\lbrack-\delta,0]$.

\noindent{\normalsize We consider the MF-DBSDE, for $r\in\left[
-\delta,0\right]  :$ }
\begin{equation}
\left\{
\begin{array}
[c]{ll}%
Y(t) & =\xi+%
%TCIMACRO{\tint _{t}^{T}}%
%BeginExpansion
{\textstyle\int_{t}^{T}}
%EndExpansion
f(s,Y_{s},Z_{s},K_{s}(\cdot),P_{(Y_{s},Z_{s},K_{s}(\cdot))})ds\\
\text{\ \ \ } & -%
%TCIMACRO{\tint _{t}^{T}}%
%BeginExpansion
{\textstyle\int_{t}^{T}}
%EndExpansion
Z(s)dB(s)-%
%TCIMACRO{\tint _{t}^{T}}%
%BeginExpansion
{\textstyle\int_{t}^{T}}
%EndExpansion%
%TCIMACRO{\tint _{E}}%
%BeginExpansion
{\textstyle\int_{E}}
%EndExpansion
K(s,\zeta)\tilde{N}(ds,d\zeta),t\in\left[  0,T\right]  ,\\
Y(t) & =Y(0),Z(t)=0,K(t,\cdot)=0,t<0,
\end{array}
\right.  \label{eq:BSFDE}%
\end{equation}
where
\[%
\begin{array}
[c]{l}%
{Y}_{s}(r,\omega):=%
\begin{cases}
Y(s+r,\omega), & r\in\lbrack-\delta,0],s+r\geq0,\omega\in\Omega,\\
Y(0,\omega), & s+r<0,
\end{cases}
\\
{Z}_{s}(r,\omega):=%
\begin{cases}
Z(s+r,\omega), & r\in\lbrack-\delta,0],s+r\geq0,\omega\in\Omega,\\
0, & s+r<0,
\end{cases}
\\
{K}_{s}(r,\omega)(\zeta):=%
\begin{cases}
K(s+r,\omega,\zeta), & r\in\lbrack-\delta,0],s+r\geq0,\omega\in\Omega,\zeta\in
E,\\
0, & s+r<0,
\end{cases}
\end{array}
\]
where $E=\mathbb{R}_{0}:=\mathbb{R}-\{0\}$ and
\[
\xi\in L^{2}(\Omega,\mathcal{F}_{T}).
\]
Here, for each $t$, $(Y_{t},Z_{t},K_{t}(\cdot))$ is assumed to belong to the
space
\[
\mathbb{S}_{\infty}^{2}\times\mathbb{L}^{2}\times\mathbb{H}^{2},
\]
of functionals defined below.

\begin{itemize}
\item $\mathbb{S}_{\infty}^{2}=L^{\infty}(\Omega,\mathcal{D}[-\delta
,0],\mathbb{R})$ consists of the càdlàg functions
\[
\alpha:[-\delta,0]\rightarrow\mathbb{R}.
\]
Let $\mathbb{S}_{\infty}^{2}$ be equipped with the norm
\[
\parallel\alpha\parallel_{\mathbb{S}_{\infty}^{2}}^{2}:=\mathbb{E}[\sup
_{r\in\lbrack-\delta,0]}|\alpha(r)|^{2}]<\infty.
\]

\end{itemize}

\noindent We refer to \cite{SEM} for more on this space in connection with
\emph{stochastic functional differential equations}.

\begin{itemize}
\item $\mathbb{L}^{2}=L_{-\delta}^{2}(\mathbb{R})$ is the space of all
functions
\[
\sigma:[-\delta,0]\rightarrow\mathbb{R},
\]
Borel measurable, such that
\[
\parallel\sigma\parallel_{\mathbb{L}^{2}}^{2}:=%
%TCIMACRO{\tint _{-\delta}^{0}}%
%BeginExpansion
{\textstyle\int_{-\delta}^{0}}
%EndExpansion
\left\vert \sigma(r)\right\vert ^{2}dr<\infty.
\]

\item $L^{2}(\nu)$ consists of Borelian functions $K:E\rightarrow\mathbb{R},$
such that
\[
\parallel K\parallel_{L^{2}(\nu)}^{2}:=%
%TCIMACRO{\tint _{E}}%
%BeginExpansion
{\textstyle\int_{E}}
%EndExpansion
K(\zeta)^{2}\nu(d\zeta)<\infty.
\]

\item $\mathbb{H}^{2}=L_{-\delta\times\nu}^{2}(\mathbb{R})$ is the space of
all functions
\[
\theta:[-\delta,0]\times E\rightarrow\mathbb{R},
\]
Borel measurable, such that
\[
\parallel\theta\parallel_{\mathbb{H}^{2}}^{2}:=%
%TCIMACRO{\tint _{-\delta}^{0}}%
%BeginExpansion
{\textstyle\int_{-\delta}^{0}}
%EndExpansion%
%TCIMACRO{\tint _{E}}%
%BeginExpansion
{\textstyle\int_{E}}
%EndExpansion
|\theta(r,\zeta)|^{2}\nu(d\zeta)dr<\infty.
\]

\item $L^{2}(\Omega,\mathcal{F}_{T})$ is the set of square integrable random
variables which are $\mathcal{F}_{T}$-measurable.
\end{itemize}

We also define the following spaces:

\begin{itemize}
\item $\mathbf{S}_{T}^{2}$ consists of the $\mathbb{F}$-adapted càdlàg
processes
\[
Y:\Omega\times\lbrack0,T]\rightarrow\mathbb{R},
\]
such that $\mathbb{E}[\sup_{t\in\lbrack0,T]}|Y(t)|^{2}]<\infty$. We equip
$\mathbf{S}_{T}^{2}$ with the norm
\[
\parallel Y\parallel_{\mathbf{S}_{T}^{2}}^{2}:=\mathbb{E}[\sup_{t\in
\lbrack0,T]}e^{\beta t}|Y(t)|^{2}],\text{ }\beta>0.
\]

\item $\mathbf{L}_{T}^{2}$ consists of the $\mathbb{F}$-predictable processes
\[
Z:\Omega\times\lbrack0,T]\rightarrow\mathbb{R},
\]
such that $\mathbb{E}[%
%TCIMACRO{\tint _{0}^{T}}%
%BeginExpansion
{\textstyle\int_{0}^{T}}
%EndExpansion
\left\vert Z(t)\right\vert ^{2}dt]<\infty$. We equip $\mathbf{L}_{T}^{2}$ with
the norm
\[
\parallel Z\parallel_{\mathbf{L}_{T}^{2}}^{2}:=\mathbb{E}[%
%TCIMACRO{\tint _{0}^{T}}%
%BeginExpansion
{\textstyle\int_{0}^{T}}
%EndExpansion
e^{\beta t}\left\vert Z(t)\right\vert ^{2}dt],\text{ }\beta>0.
\]

\item $\mathbf{H}_{T}^{2}$ consists of the $\mathbb{F}$-predictable processes
\[
K:\Omega\times\lbrack0,T]\times E\rightarrow\mathbb{R},
\]
such that $\mathbb{E}[%
%TCIMACRO{\tint _{0}^{T}}%
%BeginExpansion
{\textstyle\int_{0}^{T}}
%EndExpansion%
%TCIMACRO{\tint _{E}}%
%BeginExpansion
{\textstyle\int_{E}}
%EndExpansion
K(t,\zeta)^{2}\nu(d\zeta)dt]<\infty.$ We equip $\mathbf{H}_{T}^{2}$ with the
norm
\[
\parallel K\parallel_{\mathbf{H}_{T}^{2}}^{2}:=\mathbb{E}[%
%TCIMACRO{\tint _{0}^{T}}%
%BeginExpansion
{\textstyle\int_{0}^{T}}
%EndExpansion%
%TCIMACRO{\tint _{E}}%
%BeginExpansion
{\textstyle\int_{E}}
%EndExpansion
e^{\beta t}K(t,\zeta)^{2}\nu(d\zeta)dt],\text{ }\beta>0.
\]

\end{itemize}

\noindent Notice that if $(Y,Z,K)\in\mathbf{S}_{T,\beta}^{2}\times
\mathbf{L}_{T,\beta}^{2}\times\mathbf{H}_{T,\beta}^{2}$, then for a.e.
$t\in\lbrack-\delta,T]$, the segment process $(Y_{t},Z_{t},K_{t}(\cdot))$
belongs to $\mathbb{L}^{2}\times\mathbb{L}^{2}\times\mathbb{H}^{2}$ for a.e.
$t$, $P$-a.s.

\noindent In this section, we, as in Agram and Øksendal \cite{AO3},
\cite{AO2}, construct an Hilbert space $\mathcal{M(\cdot)}$ of (random) measures.

\begin{definition}
\noindent Let $\mathcal{\tilde{M}}(\mathcal{\mathbb{R}})$ denote the set of
random measures $\mu$ on $\mathbb{R}$ such that%
\begin{equation}
\mathbb{E[}%
%TCIMACRO{\tint _{\mathbb{R}}}%
%BeginExpansion
{\textstyle\int_{\mathbb{R}}}
%EndExpansion
|\hat{\mu}(y)|^{2}e^{-y^{2}}dy]<\infty,\label{*}%
\end{equation}
where
\[
\hat{\mu}(y)={{%
%TCIMACRO{\tint _{\mathbb{R}}}%
%BeginExpansion
{\textstyle\int_{\mathbb{R}}}
%EndExpansion
}}e^{ixy}d\mu(x)
\]
is the Fourier transform of the measure $\mu$.

\noindent If $\mu,\eta\in\mathcal{\tilde{M}}(\mathcal{\mathbb{R}})$ we define
the inner product $\left\langle \mu,\eta\right\rangle _{\mathcal{\tilde{M}%
}(\mathcal{\mathbb{R}})},$ by%
\[
\left\langle \mu,\eta\right\rangle _{\mathcal{\tilde{M}}(\mathcal{\mathbb{R}%
})}=\mathbb{E[}%
%TCIMACRO{\tint _{\mathbb{R}}}%
%BeginExpansion
{\textstyle\int_{\mathbb{R}}}
%EndExpansion
\operatorname{Re}(\overline{\hat{\mu}}(y)\hat{\eta}(y))e^{-y^{2}}dy],
\]
where, in general, $\operatorname{Re}(z)$ denotes the real part of the complex
number $z$, and $\bar{z}$ denotes the complex conjugate of $z$.

\noindent The norm $||\cdot||_{\tilde{\mathcal{M}}(\mathcal{\mathbb{R}})}$
associated to this inner product is given by
\[
\left\Vert \mu\right\Vert _{\mathcal{\tilde{M}}(\mathcal{\mathbb{R}})}%
^{2}=\left\langle \mu,\mu\right\rangle _{\mathcal{\tilde{M}}%
(\mathcal{\mathbb{R}})}=\mathbb{E[}%
%TCIMACRO{\tint _{\mathbb{R}}}%
%BeginExpansion
{\textstyle\int_{\mathbb{R}}}
%EndExpansion
|\hat{\mu}(y)|^{2}e^{-y^{2}}dy]\text{.}%
\]

\noindent The space $\mathcal{\tilde{M}}(\mathcal{\mathbb{R}})$ equipped with
the inner product $\left\langle \mu,\eta\right\rangle _{\mathcal{\tilde{M}%
}(\mathcal{\mathbb{R}})}$ is a pre-Hilbert space.

\noindent Let $\mathcal{\tilde{M}}(\mathcal{\mathbb{R}}^{m})$ denote the set
of random measures $\mu=\mu(\omega)$ on $\mathcal{\mathbb{R}}^{m}$ such that
\end{definition}

\begin{itemize}
\item We denote by $\mathcal{M}(\mathcal{\mathbb{\cdot}})$ the completion of
$\mathcal{\tilde{M}}(\mathcal{\mathbb{\cdot}})$.

\item We denote by $\mathcal{M}_{0}(\mathcal{\mathbb{\cdot}})$ the set of all
deterministic elements of $\mathcal{M}(\mathcal{\mathbb{\cdot}})$.
\end{itemize}

\begin{definition}
\begin{itemize}
\item $\mathcal{M(\mathbb{R}}^{2}\times L^{2}(\nu)\mathcal{)}$ is the space of
random measures $\mu$ on $\mathcal{\mathbb{R}}^{2}\times L^{2}(\nu),$ such
that
\[%
\begin{array}
[c]{lll}%
\left\Vert \mu\right\Vert _{\mathcal{M(\mathcal{\mathbb{R}}}^{2}\times
L^{2}(\nu)\mathcal{)}}^{2} & := & \mathbb{E[}%
%TCIMACRO{\tint _{\mathcal{\mathbb{R}}_{0}}}%
%BeginExpansion
{\textstyle\int_{\mathcal{\mathbb{R}}_{0}}}
%EndExpansion%
%TCIMACRO{\tint _{\mathcal{\mathbb{R}}^{3}}}%
%BeginExpansion
{\textstyle\int_{\mathcal{\mathbb{R}}^{3}}}
%EndExpansion
|\widehat{\mu}(y_{1},y_{2},y_{3},\zeta)|^{2}\exp(-\underset{j=1}{\overset{3}{%
%TCIMACRO{\tsum }%
%BeginExpansion
{\textstyle\sum}
%EndExpansion
}}y_{j}^{2})\nu d(y_{1},y_{2},y_{3})(d\zeta)]<\infty\text{,}%
\end{array}
\]
where $\widehat{\mu}(y_{1},y_{2},y_{3},\zeta)=\widehat{\mu}(y_{1},y_{2}%
,y_{3})$ is the Fourier transform of the measure $\mu$ parametrized at $\zeta
$, i.e.,%
\[%
\begin{array}
[c]{lll}%
\widehat{\mu}(y_{1},y_{2},y_{3}) & := & {%
%TCIMACRO{\tint _{\mathcal{\mathbb{R}}^{3}}}%
%BeginExpansion
{\textstyle\int_{\mathcal{\mathbb{R}}^{3}}}
%EndExpansion
}\exp(-2\pi i(\underset{j=1}{\overset{3}{%
%TCIMACRO{\tsum }%
%BeginExpansion
{\textstyle\sum}
%EndExpansion
}}x_{j}y_{j}))d\mu(x_{1},x_{2},x_{3});\quad y_{1},y_{2},y_{3}\in\mathbb{R}.
\end{array}
\]

\end{itemize}
\end{definition}

\begin{lemma}
Let $X^{(1)},X^{(2)},\widetilde{X}^{(1)},\widetilde{X}^{(2)}$ and
$X^{(3)},\widetilde{X}^{(3)}$ be random variables in $L^{2}(P)$ and in
$L^{2}(\nu)$ respectively. Then
\begin{align*}
&  ||\mathcal{L}(X^{(1)},X^{(2)},X^{(3)})-\mathcal{L}(\widetilde{X}%
^{(1)},\widetilde{X}^{(2)},\widetilde{X}^{(3)})||_{\mathcal{M}_{0}%
\mathcal{(\mathbb{R}}^{2}\times L^{2}(\nu)\mathcal{)}}^{2}\\
&  \leq C\mathbb{E}[(X^{(1)}-\widetilde{X}^{(1)})^{2}+(X^{(2)},\widetilde
{X}^{(2)})^{2}+%
%TCIMACRO{\tint _{\mathbb{R}_{0}}}%
%BeginExpansion
{\textstyle\int_{\mathbb{R}_{0}}}
%EndExpansion
(X^{(3)}(\zeta),\widetilde{X}^{(3)}(\zeta))^{2}\nu(\zeta)]\text{,}%
\end{align*}
where $\mathcal{L}(X)=P_{X}.$
\end{lemma}

\noindent{Proof}\quad Let $X=(X^{(1)},X^{(2)},X^{(3)}),$ $\widetilde
{X}=(\widetilde{X}^{(1)},\widetilde{X}^{(2)},\widetilde{X}^{(3)})$ and
$y=(y_{1},y_{2},y_{3}):$%
\begin{align*}
&  ||\mathcal{L}(X)-\mathcal{L}(\widetilde{X})||_{\mathcal{M}_{0}%
\mathcal{(\mathbb{R}}^{2}\times L^{2}(\nu)\mathcal{)}}^{2}\\
&  =%
%TCIMACRO{\tint _{\mathcal{\mathbb{R}}_{0}}}%
%BeginExpansion
{\textstyle\int_{\mathcal{\mathbb{R}}_{0}}}
%EndExpansion%
%TCIMACRO{\tint _{\mathcal{\mathbb{R}}^{3}}}%
%BeginExpansion
{\textstyle\int_{\mathcal{\mathbb{R}}^{3}}}
%EndExpansion
|\widehat{\mathcal{L}}(X)(y,\zeta)-\widehat{\mathcal{L}}(\widetilde
{X})(y,\zeta)|^{2}e^{-y^{2}}dy\min(1,\zeta^{2})\nu(\zeta)\\
&  =%
%TCIMACRO{\tint _{\mathcal{\mathbb{R}}_{0}}}%
%BeginExpansion
{\textstyle\int_{\mathcal{\mathbb{R}}_{0}}}
%EndExpansion%
%TCIMACRO{\tint _{\mathcal{\mathbb{R}}^{3}}}%
%BeginExpansion
{\textstyle\int_{\mathcal{\mathbb{R}}^{3}}}
%EndExpansion
|%
%TCIMACRO{\tint _{\mathcal{\mathbb{R}}^{3}}}%
%BeginExpansion
{\textstyle\int_{\mathcal{\mathbb{R}}^{3}}}
%EndExpansion
e^{-ixy}d\mathcal{L}(X)(x)-%
%TCIMACRO{\tint _{\mathcal{\mathbb{R}}^{3}}}%
%BeginExpansion
{\textstyle\int_{\mathcal{\mathbb{R}}^{3}}}
%EndExpansion
e^{-ixy}d\mathcal{L}(\widetilde{X})(x)|^{2}e^{-y^{2}}dy\min(1,\zeta^{2}%
)\nu(\zeta)\\
&  \leq%
%TCIMACRO{\tint _{\mathcal{\mathbb{R}}_{0}}}%
%BeginExpansion
{\textstyle\int_{\mathcal{\mathbb{R}}_{0}}}
%EndExpansion%
%TCIMACRO{\tint _{\mathcal{\mathbb{R}}^{3}}}%
%BeginExpansion
{\textstyle\int_{\mathcal{\mathbb{R}}^{3}}}
%EndExpansion
|\mathbb{E}[e^{-iXy}-e^{-i\widetilde{X}y}]|^{2}e^{-y^{2}}dy\min(1,\zeta
^{2})\nu(\zeta)\\
&  \leq%
%TCIMACRO{\tint _{\mathcal{\mathbb{R}}_{0}}}%
%BeginExpansion
{\textstyle\int_{\mathcal{\mathbb{R}}_{0}}}
%EndExpansion%
%TCIMACRO{\tint _{\mathcal{\mathbb{R}}^{3}}}%
%BeginExpansion
{\textstyle\int_{\mathcal{\mathbb{R}}^{3}}}
%EndExpansion
|\mathbb{E}[(X-\widetilde{X})y]|^{2}e^{-y^{2}}dy\min(1,\zeta^{2})\nu(\zeta)\\
&  \leq%
%TCIMACRO{\tint _{\mathcal{\mathbb{R}}_{0}}}%
%BeginExpansion
{\textstyle\int_{\mathcal{\mathbb{R}}_{0}}}
%EndExpansion%
%TCIMACRO{\tint _{\mathcal{\mathbb{R}}^{3}}}%
%BeginExpansion
{\textstyle\int_{\mathcal{\mathbb{R}}^{3}}}
%EndExpansion
y^{2}\mathbb{E}[(X-\widetilde{X})^{2}]e^{-y^{2}}dy\min(1,\zeta^{2})\nu
(\zeta)\\
&  =\mathbb{E}[(X^{(1)}-\widetilde{X}^{(1)})^{2}]%
%TCIMACRO{\tint _{\mathcal{\mathbb{R}}_{0}}}%
%BeginExpansion
{\textstyle\int_{\mathcal{\mathbb{R}}_{0}}}
%EndExpansion%
%TCIMACRO{\tint _{\mathcal{\mathbb{R}}^{3}}}%
%BeginExpansion
{\textstyle\int_{\mathcal{\mathbb{R}}^{3}}}
%EndExpansion
y^{2}e^{-y^{2}}dy\min(1,\zeta^{2})\nu(\zeta)\\
&  +\mathbb{E}[(X^{(2)},\widetilde{X}^{(2)})^{2}]%
%TCIMACRO{\tint _{\mathcal{\mathbb{R}}_{0}}}%
%BeginExpansion
{\textstyle\int_{\mathcal{\mathbb{R}}_{0}}}
%EndExpansion%
%TCIMACRO{\tint _{\mathcal{\mathbb{R}}^{3}}}%
%BeginExpansion
{\textstyle\int_{\mathcal{\mathbb{R}}^{3}}}
%EndExpansion
y^{2}e^{-y^{2}}dy\min(1,\zeta^{2})\nu(\zeta)\\
&  +%
%TCIMACRO{\tint _{\mathbb{R}_{0}}}%
%BeginExpansion
{\textstyle\int_{\mathbb{R}_{0}}}
%EndExpansion
\mathbb{E}[(X^{(3)},\widetilde{X}^{(3)})^{2}(\zeta)]\min(1,\zeta^{2})\nu
(\zeta)%
%TCIMACRO{\tint _{\mathcal{\mathbb{R}}^{3}}}%
%BeginExpansion
{\textstyle\int_{\mathcal{\mathbb{R}}^{3}}}
%EndExpansion
y^{2}e^{-y^{2}}dy.
\end{align*}
Note that $%
%TCIMACRO{\tint _{\mathcal{\mathbb{R}}_{0}}}%
%BeginExpansion
{\textstyle\int_{\mathcal{\mathbb{R}}_{0}}}
%EndExpansion
\min(1,\zeta^{2})\nu(\zeta)<\infty$ for all Lévy measure $\nu.$\textbf{$\qquad
\qquad$}$\square$

\begin{definition}
\begin{itemize}
\item Define $\mathcal{M}^{\delta}\mathcal{(\mathbb{R}}^{2}\times L^{2}%
(\nu)\mathcal{)}$ to be the Hilbert space of all path segments $\overline{\mu
}=\{\mu(s)\}_{s\in\ [0,\delta]}$ of measure-valued processes $\mu(\cdot)$ with
$\mu(s)\in\mathcal{M(\mathbb{R}}^{2}\times L^{2}(\nu)\mathcal{)}$ for each
$s\in\lbrack-\delta,0]$, equipped with the norm
\[
\left\Vert \overline{\mu}\right\Vert _{\mathcal{M}^{\delta}%
\mathcal{(\mathbb{R}}^{2}\times L^{2}(\nu)\mathcal{)}}:={%
%TCIMACRO{\tint _{-\delta}^{0}}%
%BeginExpansion
{\textstyle\int_{-\delta}^{0}}
%EndExpansion
}\left\Vert \mu(s)\right\Vert _{\mathcal{M(\mathbb{R}}^{2}\times L^{2}%
(\nu)\mathcal{)}}ds.
\]

\item $\mathcal{M}_{0}\mathcal{(\mathbb{R}}^{2}\times L^{2}(\nu)\mathcal{)}$
and $\mathcal{M}_{0}^{\delta}\mathcal{(\mathbb{R}}^{2}\times L^{2}%
(\nu)\mathcal{)}$ denote the set of deterministic elements of
$\mathcal{M(\mathbb{R}}^{2}\times L^{2}(\nu)\mathcal{)}$ and $\mathcal{M}%
^{\delta}\mathcal{(\mathbb{R}}^{2}\times L^{2}(\nu)\mathcal{)}$,
respectively.\newline
\end{itemize}
\end{definition}

\noindent There are several advantages with working with this Hilbert space
$\mathcal{M}$, compared to the Wasserstein metric space:

\begin{itemize}
\item A Hilbert space has a useful stronger structure than a metric space.

\item The distance is not continuous but the norm is.

\item The Wasserstein metric space $\mathcal{P}_{2}$ deals only with
probability measures with finite second moment, while the Hilbert space deals
with any (random) measure satisfying (\ref{*}).
\end{itemize}

\noindent Let us give some examples of measures:

\begin{example}
[Measures]\-We consider here a $1$-dimentional case:

\begin{enumerate}
\item Suppose that $\mu=\delta_{x_{0}}$, the unit point mass at $x_{0}%
\in\mathbb{R}$. Then $\delta_{x_{0}}\in\mathcal{M}_{0}(\mathbb{R})$ and
\[
{%
%TCIMACRO{\tint _{\mathbb{R}}}%
%BeginExpansion
{\textstyle\int_{\mathbb{R}}}
%EndExpansion
}e^{ixy}d\mu(x)=e^{ix_{0}y},
\]

and hence
\[%
\begin{array}
[c]{lll}%
\left\Vert \mu\right\Vert _{\mathcal{M}_{0}(\mathbb{R})}^{2} & =%
%TCIMACRO{\tint _{\mathbb{R}}}%
%BeginExpansion
{\textstyle\int_{\mathbb{R}}}
%EndExpansion
|e^{ix_{0}y}|^{2}e^{-y^{2}}dy & <\infty\text{.}%
\end{array}
\]

\item Suppose $d\mu(x)=f(x)dx$, where $f\in L^{1}(\mathbb{R})$. Then $\mu
\in\mathcal{M}_{0}(\mathbb{R})$ and by Riemann-Lebesque lemma, $\hat{\mu
}(y)\in C_{0}(\mathbb{R})$, i.e. $\hat{\mu}$ is continuous and $\hat{\mu
}(y)\rightarrow0$ when $|y|\rightarrow\infty$. In particular, $|\hat{\mu}|$ is
bounded on $\mathbb{R}$ and hence%
\[%
\begin{array}
[c]{lll}%
\left\Vert \mu\right\Vert _{\mathcal{M}_{0}(\mathbb{R})}^{2} & =%
%TCIMACRO{\tint _{\mathbb{R}}}%
%BeginExpansion
{\textstyle\int_{\mathbb{R}}}
%EndExpansion
|\hat{\mu}(y)|^{2}e^{-y^{2}}dy & <\infty\text{.}%
\end{array}
\]

\item Suppose that $\mu$ is any finite positive measure on $\mathbb{R}$. Then
$\mu\in\mathcal{M}_{0}(\mathbb{R})$ and
\[%
\begin{array}
[c]{lll}%
|\hat{\mu}(y)| & \leq%
%TCIMACRO{\tint _{\mathbb{R}}}%
%BeginExpansion
{\textstyle\int_{\mathbb{R}}}
%EndExpansion
d\mu(y)=\mu(\mathbb{R}) & <\infty\text{, for all }y\text{,}%
\end{array}
\]
and hence%
\[%
\begin{array}
[c]{lll}%
\left\Vert \mu\right\Vert _{\mathcal{M}_{0}(\mathbb{R})}^{2} & =%
%TCIMACRO{\tint _{\mathbb{R}}}%
%BeginExpansion
{\textstyle\int_{\mathbb{R}}}
%EndExpansion
|\hat{\mu}(y)|^{2}e^{-y^{2}}dy & <\infty\text{.}%
\end{array}
\]

\item Next, suppose $x_{0}=x_{0}(\omega)$ is random. Then $\delta
_{x_{0}(\omega)}$ is a random measure in $\mathcal{M}(\mathbb{R})$. Similarly,
if $f(x)=f(x,\omega)$ is random, then $d\mu(x,\omega)=f(x,\omega)dx$ is a
random measure in $\mathcal{M}(\mathbb{R})$.\newline
\end{enumerate}
\end{example}

\section{Existence and uniqueness}

\noindent The objective of this section is to give a theorem on the existence
and the uniqueness of the solution of equation (\ref{eq:BSFDE}). For the proof
we use the fixed point argument which is a classical tool to prove the
existence and the uniqueness of BSDEs.

\begin{definition}
\-

\begin{itemize}
\item A functional
\[
f:\Omega\mathcal{\times}[0,T]\times\mathbb{L}^{2}\times\mathbb{L}^{2}%
\times\mathbb{H}^{2}\times\mathcal{M}^{\delta}\mathcal{(\mathbb{R}}^{2}\times
L^{2}(\nu)\mathcal{)}\rightarrow\mathbb{R},
\]
is progressively measurable.

\item A process
\[
(Y,Z,K)\in\mathbf{S}_{T}^{2}\times\mathbf{L}_{T}^{2}\times\mathbf{H}_{T}^{2}%
\]
is said to be a \emph{solution} to (\ref{eq:BSFDE}) if%
\[%
%TCIMACRO{\tint _{0}^{T}}%
%BeginExpansion
{\textstyle\int_{0}^{T}}
%EndExpansion
\left\vert f(s,{Y}_{s},{Z}_{s},{K}_{s}(\cdot),P_{(Y_{s},Z_{s},K_{s}(\cdot
))})\right\vert ds<+\infty\text{ }P\text{-a.s.,}%
\]
and%
\[
\left\{
\begin{array}
[c]{ll}%
Y(t) & =\xi+%
%TCIMACRO{\tint _{t}^{T}}%
%BeginExpansion
{\textstyle\int_{t}^{T}}
%EndExpansion
f(s,{Y}_{s},{Z}_{s},{K}_{s}(\cdot),P_{(Y_{s},Z_{s},K_{s}(\cdot))})ds-%
%TCIMACRO{\tint _{t}^{T}}%
%BeginExpansion
{\textstyle\int_{t}^{T}}
%EndExpansion
Z(s)dB(s)\\
& \text{ \ \ \ \ \ }-%
%TCIMACRO{\tint _{t}^{T}}%
%BeginExpansion
{\textstyle\int_{t}^{T}}
%EndExpansion%
%TCIMACRO{\tint _{E}}%
%BeginExpansion
{\textstyle\int_{E}}
%EndExpansion
K(s,\zeta)\tilde{N}(ds,d\zeta),t\in\lbrack0,T],\\
Y(t) & =Y(0),Z(0)=K(t,\cdot)=0,t<0.
\end{array}
\right.
\]
We impose the following set of assumptions which will garantee the existence
and the uniqueness of the solution of the MF-DBSDE (\ref{eq:BSFDE}).
\end{itemize}
\end{definition}

\noindent{\normalsize \textbf{Assumptions }}

\noindent Let $f$ be a functional generator and $\xi$ the terminal condition.
Suppose that:

\begin{description}
\item[(i)] $\xi\in L^{2}\left(  \Omega,\mathcal{F}_{T}\right)  $.

\item[(ii)] For all $t\in\lbrack0,T],$ we have%
\[
|f(t,0,0,0,P_{0})|<c,
\]
where $P_{0}$ is the Dirac measure with mass at zero and $c$ is a given constant.

\item[(iii)] {\normalsize For all $t\in\left[  0,T\right]  $ and for all
$y_{i}\in\mathbb{L}^{2},z_{i}\in\mathbb{L}^{2},k_{i}$}$(\cdot)\in
\mathbb{H}^{2}$ {\normalsize and }$\eta_{i}\in\mathcal{M}_{0}(\mathbb{L}%
^{2}\times\mathbb{L}^{2}\times\mathbb{H}^{2}\mathcal{)},i=1,2${\normalsize $,$%
\ we have for a constant $C>0$ and for some probability measure }$\mu$
{\normalsize on }$[-\delta,0]\times\mathcal{B}[-\delta,0]$ where $\mathcal{B}$
stands for the Borel sets of $[-\delta,0]$, such that
\begin{align*}
&  \left\vert f(t,y_{1},z_{1},k_{1}(\cdot),\eta_{1})-f(t,y_{2},z_{2}%
,k_{2}(\cdot),\eta_{2})\right\vert ^{2}\\
&  \leq C%
%TCIMACRO{\tint _{-\delta}^{0}}%
%BeginExpansion
{\textstyle\int_{-\delta}^{0}}
%EndExpansion
(|Y_{1}(t+r)-Y_{2}(t+r)|^{2}+|Z_{1}(t+r)-Z_{2}(t+r)|^{2}\\
&  +%
%TCIMACRO{\tint _{E}}%
%BeginExpansion
{\textstyle\int_{E}}
%EndExpansion
|K_{1}(t+r,\zeta)-K_{2}(t+r,\zeta)|^{2}\nu(d\zeta)\\
&  +||\eta_{1}(r)-\eta_{2}(r)||_{\mathcal{M(\mathbb{R}}^{2}\times L^{2}%
(\nu)\mathcal{)}}^{2})\mu(dr),\text{ }P\text{-a.s.}%
\end{align*}

\end{description}

\noindent The following theorem gives the existence and the uniqueness of the
solution of a MF-DBSDE with jumps under assumptions (i)-(iii).

\begin{theorem}
\label{main} Let us suppose the above assumptions (i)-(iii), with $\rho
>\mu(\{0\}).$ Then for sufficiently small $\delta_{\rho}>0,$ it holds for all
$\delta\in\left(  0,\delta_{\rho}\right)  ,$ the MF-DBSDE \ (\ref{eq:BSFDE})
admits a unique solution\textbf{\ }$(Y,Z,K)\in\mathbf{S}_{T}^{2}%
\times\mathbf{L}_{T}^{2}\times\mathbf{H}_{T}^{2}\mathbf{.}$
\end{theorem}

\begin{remark}
In general, the Lipschitz condition with $||y_{1}-y_{2}||_{\mathbb{S}_{\infty
}^{2}}^{2}$ instead of\newline$||y_{1}-y_{2}||_{\mathbb{L}^{2}}^{2}$ is
considered; see, for instance Delong and Imkeller \cite{di}, \cite{DI2}.
Recall that:
\[
||y_{1}-y_{2}||_{\mathbb{L}^{2}}^{2}\leq C||y_{1}-y_{2}||_{\mathbb{S}_{\infty
}^{2}}^{2},
\]
they obtain the solution of the DBSDE for a sufficiently small time horizon
$T>0$ or a sufficiently small Lipschitz constants. We take here a condition
which is more restrictive but it allows us to obtain the existence and
uniqueness of our MF-DBSDE with jumps (\ref{eq:BSFDE}) for any finite time
horizon $T>0$ and for any Lipschitz constant but for a sufficiently small
delay constant $\delta>0$.
\end{remark}

\noindent{Proof}\quad Let us define the mapping\textbf{\ }%
\[
\Phi:\left(  L^{2}(\mathcal{F}_{0})\times\mathbf{L}_{T}^{2}\right)
\times\mathbf{L}_{T}^{2}\times\mathbf{H}_{T}^{2}\rightarrow\left(
L^{2}(\mathcal{F}_{0})\times\mathbf{L}_{T}^{2}\right)  \times\mathbf{L}%
_{T}^{2}\times\mathbf{H}_{T}^{2}%
\]
by setting\
\[
\Phi\left(  \left(  U(0),U\right)  ,V,Q\right)  :=\left(  \left(
Y(0),Y\right)  ,Z,K\right)  ,
\]
where for\ $U(0)\in L^{2}(\mathcal{F}_{0}),$ $U\in\mathbf{L}_{T}^{2},$
$U(t)=U(0),$ $t<0,$ $V\in\mathbf{L}_{T}^{2},$ $V(t)=0,$ $t<0,$ $Q\in
\mathbf{H}_{T}^{2},$ $Q(t)=0,$ $t<0,$ $(Y,Z,K)\in\mathbf{S}_{T}^{2}%
\times\mathbf{L}_{T}^{2}\times\mathbf{H}_{T}^{2}$ $(\subset(L^{2}%
(\mathcal{F}_{0})\times\mathbf{L}_{T}^{2})\times\mathbf{L}_{T}^{2}%
\times\mathbf{H}_{T}^{2})$ is the unique solution of the MF-DBSDE with jumps%
\[
\left\{
\begin{array}
[c]{ll}%
Y(t) & =\xi+\int_{t}^{T}f(s,U_{s},V_{s},Q_{s}(\cdot),P_{(U_{s},V_{s}%
,Q_{s}(\cdot))})ds-\int_{t}^{T}Z(s)dB(s)\\
& \text{ \ \ \ \ \ \ \ \ \ \ \ \ }-\int_{t}^{T}\int_{E}K(s,\zeta)\tilde
{N}(ds,d\zeta),t\in\left[  0,T\right]  ,\\
Y(t) & =Y(0),\text{ }Z(t)=0,\text{ }K(t,\cdot)=0,t<0.
\end{array}
\right.
\]
For\ $\beta>0$ and $\left(  \left(  U(0),U\right)  ,V,Q\right)  \in
((L^{2}(\mathcal{F}_{0})\times\mathbf{L}_{T}^{2})\times\mathbf{L}_{T}%
^{2}\times\mathbf{H}_{T}^{2})$ we introduce the norm%
\[%
\begin{array}
[c]{ll}%
\left\Vert \left(  \left(  U(0),U\right)  ,V,Q\right)  \right\Vert  &
:=\left\Vert \left(  \left(  U(0),U\right)  ,V,Q\right)  \right\Vert _{\beta
}\\
& :=(\mathbb{E}[\left\vert U(0)\right\vert ^{2}]\\
& +\mathbb{E[}\int_{0}^{T}e^{\beta s}(\left\vert U(s)\right\vert
^{2}+\left\vert V(s)\right\vert ^{2}+\int_{E}\left\vert Q(s,\zeta)\right\vert
^{2}\nu(d\zeta))ds])^{\frac{1}{2}}.
\end{array}
\]
Note that $(L^{2}(\mathcal{F}_{0})\times\mathbf{L}_{T}^{2})\times
\mathbf{L}_{T}^{2}\times\mathbf{H}_{T}^{2}$ endowed with this norm is a Banach
space. We will show that for suitably chosen\ $\beta>0,$ $\delta\in\left(
0,\delta_{0}\right)  $\textbf{$,$ }the mapping\textbf{\ }%
\[
\Phi:(\left(  L^{2}(\mathcal{F}_{0})\times\mathbf{L}_{T}^{2}\right)
\times\mathbf{L}_{T}^{2}\times\mathbf{H}_{T}^{2},\left\Vert \cdot\right\Vert
_{\beta})\mathbf{\rightarrow}\left(  L^{2}(\mathcal{F}_{0})\times
\mathbf{L}_{T}^{2}\right)  \times\mathbf{L}_{T}^{2}\times\mathbf{H}_{T}%
^{2},\left\Vert \cdot\right\Vert _{\beta})
\]
\ is contracting, i.e., there is a unique fixed point\ $\left(  \left(
Y(0),Y\right)  ,Z,K\right)  \in$\newline$\left(  \left(  L^{2}(\mathcal{F}%
_{0})\times\mathbf{L}_{T}^{2}\right)  \text{ }\times\mathbf{L}_{T}^{2}%
\times\mathbf{H}_{T}^{2}\right)  $ of\ $\Phi.$ Consequently,\textbf{\ }%
\[
\left\{
\begin{array}
[c]{ll}%
Y(t) & =\xi+\int_{t}^{T}f(s,Y_{s},Z_{s},K_{s}(\cdot),P_{(Y_{s},Z_{s}%
,K_{s}(\cdot))})ds-\int_{t}^{T}Z(s)dB(s)\\
& \text{ \ \ \ \ \ \ \ \ \ \ \ \ \ }-\int_{t}^{T}\int_{E}K(s,\zeta)\tilde
{N}(ds,d\zeta),t\in\left[  0,T\right]  ,\\
Y(t) & =Y(0),\text{ }Z(t)=0,\text{ }K(t,\cdot)=0,t<0.
\end{array}
\right.
\]
In particular $Y$ has a continuous version and by standard estimations, there
exists a constant $C\mathbf{\in%
%TCIMACRO{\U{211d} }%
%BeginExpansion
\mathbb{R}
%EndExpansion
,}$ such that%
\[%
\begin{array}
[c]{l}%
\mathbb{E}[\underset{t\in\left[  0,T\right]  }{\sup}\left\vert Y(t)\right\vert
^{2}]\leq C(\mathbb{E}[\left\vert \xi\right\vert ^{2}]+\mathbb{E}[\int_{t}%
^{T}\left\vert f(s,U_{s},V_{s},Q_{s}(\cdot),P_{(U_{s},V_{s},Q_{s}(\cdot
))})\right\vert ^{2}ds]\\
\text{ \ \ \ \ \ \ \ \ \ \ \ \ \ \ \ \ \ \ \ \ \ \ \ \ }+\mathbb{E[}\int
_{t}^{T}\left\vert Z(s)\right\vert ^{2}ds]+\mathbb{E[}\int_{t}^{T}\int
_{E}\left\vert K(s,\zeta)\right\vert ^{2}\nu(d\zeta)ds])<\infty.
\end{array}
\]
Consequently, $Y$\textbf{$\in\mathbf{S}_{T}^{2}.$ }Let us consider\textbf{\ }%
$(\left(  U(0),U\right)  ,V,Q),(\left(  U^{\prime}(0),U^{\prime}\right)
,V^{\prime},Q^{\prime})\in\left(  L^{2}(\mathcal{F}_{0})\mathbf{\times
\mathbf{L}_{T}^{2}}\right)  $\textbf{$\times\mathbf{L}_{T}^{2}\times H_{T}%
^{2}$ }and let us use the simplified notations:%
\[%
\begin{array}
[c]{ll}%
\Phi(\left(  U(0),U\right)  ,V,Q) & :=\left(  \left(  Y(0),Y\right)
,Z,K\right)  ,\\
\Phi(\left(  U^{\prime}(0),U^{\prime}\right)  ,V^{\prime},Q^{\prime}) &
:=\left(  \left(  Y^{\prime}(0),Y^{\prime}\right)  ,Z^{\prime},K^{\prime
}\right)  ,\\
(\left(  \bar{U}(0),\bar{U}\right)  ,\bar{V},\bar{Q}) & :=(\left(
U(0),U\right)  ,V,Q)-(\left(  U^{\prime}(0),U^{\prime}\right)  ,V^{\prime
},Q^{\prime}),\\
\left(  \left(  \bar{Y}(0),\bar{Y}\right)  ,\bar{Z},\bar{K}\right)  &
:=\left(  \left(  Y(0),Y\right)  ,Z,K\right)  -\left(  \left(  Y^{\prime
}(0),Y^{\prime}\right)  ,Z^{\prime},K^{\prime}\right)  ,
\end{array}
\]
Applying Itô's formula to\textbf{\ }$(e^{\beta t}\left\vert \bar
{Y}(t)\right\vert ^{2})_{t\geq0}$\textbf{\ }and using the Lipschitz condition
$(ii)$, we get%
\begin{equation}%
\begin{array}
[c]{l}%
e^{\beta t}\left\vert \bar{Y}(t)\right\vert ^{2}+%
%TCIMACRO{\tint _{t}^{T}}%
%BeginExpansion
{\textstyle\int_{t}^{T}}
%EndExpansion
e^{\beta s}(\beta\left\vert \bar{Y}(s)\right\vert ^{2}+\left\vert \bar
{Z}(s)\right\vert ^{2}+%
%TCIMACRO{\tint _{E}}%
%BeginExpansion
{\textstyle\int_{E}}
%EndExpansion
\left\vert \bar{K}(s,\zeta)\right\vert ^{2}\nu(d\zeta))ds\\
\leq2%
%TCIMACRO{\tint _{t}^{T}}%
%BeginExpansion
{\textstyle\int_{t}^{T}}
%EndExpansion
e^{\beta s}\left\vert \bar{Y}(s)\right\vert \times\\
\times|f(s,U_{s},V_{s},Q_{s}(\cdot),P_{(U_{s},V_{s},Q_{s}(\cdot))}%
)-f(s,U_{s}^{\prime},V_{s}^{\prime},Q_{s}^{\prime}(\cdot),P_{(U_{s}^{\prime
},V_{s}^{\prime},Q_{s}^{\prime}(\cdot))})|ds\\
-2%
%TCIMACRO{\tint _{t}^{T}}%
%BeginExpansion
{\textstyle\int_{t}^{T}}
%EndExpansion
e^{\beta s}\bar{Y}(s)\cdot\bar{Z}(s)dB(s)-2%
%TCIMACRO{\tint _{t}^{T}}%
%BeginExpansion
{\textstyle\int_{t}^{T}}
%EndExpansion%
%TCIMACRO{\tint _{E}}%
%BeginExpansion
{\textstyle\int_{E}}
%EndExpansion
e^{\beta s}\bar{Y}(s^{-})\cdot\bar{K}(s,\zeta)\tilde{N}(ds,d\zeta)\\
\leq2%
%TCIMACRO{\tint _{t}^{T}}%
%BeginExpansion
{\textstyle\int_{t}^{T}}
%EndExpansion
e^{\beta s}\left\vert \bar{Y}(s)\right\vert \times\\
\times C(%
%TCIMACRO{\tint _{-\delta}^{0}}%
%BeginExpansion
{\textstyle\int_{-\delta}^{0}}
%EndExpansion
\{\left\vert \bar{U}(s+r)\right\vert +\left\vert \bar{V}(s+r)\right\vert +(%
%TCIMACRO{\tint _{E}}%
%BeginExpansion
{\textstyle\int_{E}}
%EndExpansion
\left\vert \bar{Q}(s+r,\zeta)\right\vert ^{2}\nu(d\zeta))^{\frac{1}{2}}%
\}\mu(dr)\\
+\sqrt{\pi}\mathbb{E[}%
%TCIMACRO{\tint _{-\delta}^{0}}%
%BeginExpansion
{\textstyle\int_{-\delta}^{0}}
%EndExpansion
\{\left\vert \bar{U}(s+r)\right\vert +\left\vert \bar{V}(s+r)\right\vert +(%
%TCIMACRO{\tint _{E}}%
%BeginExpansion
{\textstyle\int_{E}}
%EndExpansion
\left\vert \bar{Q}(s+r,\zeta)\right\vert ^{2}\nu(d\zeta))^{\frac{1}{2}}%
\}\mu(dr)])ds\\
-2%
%TCIMACRO{\tint _{t}^{T}}%
%BeginExpansion
{\textstyle\int_{t}^{T}}
%EndExpansion
e^{\beta s}\bar{Y}(s)\cdot\bar{Z}(s)dB(s)-2%
%TCIMACRO{\tint _{t}^{T}}%
%BeginExpansion
{\textstyle\int_{t}^{T}}
%EndExpansion%
%TCIMACRO{\tint _{E}}%
%BeginExpansion
{\textstyle\int_{E}}
%EndExpansion
e^{\beta s}\bar{Y}(s^{-})\cdot\bar{K}(s,\zeta)\tilde{N}(ds,d\zeta).
\end{array}
\label{ra}%
\end{equation}
Using for every term in the integrand of the Lebesgue integral at the right
hand side of the above equality, the estimate $2C(1+%
%TCIMACRO{\tint _{\mathbb{R} ^{3}}}%
%BeginExpansion
{\textstyle\int_{\mathbb{R} ^{3}}}
%EndExpansion
y^{2}e^{-y^{2}}dy)ab\leq2\rho C^{\prime2}a^{2}+\tfrac{1}{\rho}b^{2}$, we
obtain%
\[%
\begin{array}
[c]{l}%
e^{\beta t}\left\vert \bar{Y}(t)\right\vert ^{2}+%
%TCIMACRO{\tint _{t}^{T}}%
%BeginExpansion
{\textstyle\int_{t}^{T}}
%EndExpansion
e^{\beta s}(\beta\left\vert \bar{Y}(s)\right\vert ^{2}+\left\vert \bar
{Z}(s)\right\vert ^{2}+%
%TCIMACRO{\tint _{E}}%
%BeginExpansion
{\textstyle\int_{E}}
%EndExpansion
\left\vert \bar{K}(s,\zeta)\right\vert ^{2}\nu(d\zeta))ds\\
\leq12\rho C^{\prime2}%
%TCIMACRO{\tint _{t}^{T}}%
%BeginExpansion
{\textstyle\int_{t}^{T}}
%EndExpansion
e^{\beta s}\left\vert \bar{Y}(s)\right\vert ^{2}ds\\
+\tfrac{1}{\rho}%
%TCIMACRO{\tint _{t}^{T}}%
%BeginExpansion
{\textstyle\int_{t}^{T}}
%EndExpansion
e^{\beta s}%
%TCIMACRO{\tint _{-\delta}^{0}}%
%BeginExpansion
{\textstyle\int_{-\delta}^{0}}
%EndExpansion
\{\left\vert \bar{U}(s+r)\right\vert ^{2}+\left\vert \bar{V}(s+r)\right\vert
^{2}+%
%TCIMACRO{\tint _{E}}%
%BeginExpansion
{\textstyle\int_{E}}
%EndExpansion
\left\vert \bar{Q}(s+r,\zeta)\right\vert ^{2}\nu(d\zeta)\}\mu(dr)ds\\
-2%
%TCIMACRO{\tint _{t}^{T}}%
%BeginExpansion
{\textstyle\int_{t}^{T}}
%EndExpansion
e^{\beta s}\bar{Y}(s)\cdot\bar{Z}(s)dB(s)-2%
%TCIMACRO{\tint _{t}^{T}}%
%BeginExpansion
{\textstyle\int_{t}^{T}}
%EndExpansion%
%TCIMACRO{\tint _{E}}%
%BeginExpansion
{\textstyle\int_{E}}
%EndExpansion
e^{\beta s}\bar{Y}(s)\cdot\bar{K}(s,\zeta)\tilde{N}(ds,d\zeta),
\end{array}
\]
{where }$C^{\prime2}=C^{2}(1+(%
%TCIMACRO{\tint _{\mathbb{R} ^{3}}}%
%BeginExpansion
{\textstyle\int_{\mathbb{R} ^{3}}}
%EndExpansion
y^{2}e^{-y^{2}}dy)^{2}+2%
%TCIMACRO{\tint _{\mathbb{R} ^{3}}}%
%BeginExpansion
{\textstyle\int_{\mathbb{R} ^{3}}}
%EndExpansion
y^{2}e^{-y^{2}}dy.$

\noindent Choose\textbf{\ }$\beta=1+12\rho C^{\prime2},$\textbf{\ }the last
line in (\ref{ra}) constitutes a sum of martingale differences (with the
necessary integrability properties)\textbf{\ }then, taking expectation, we get%
\begin{equation}%
\begin{array}
[c]{l}%
\mathbb{E}[e^{\beta t}\left\vert \bar{Y}(t)\right\vert ^{2}]+\mathbb{E}[%
%TCIMACRO{\tint _{t}^{T}}%
%BeginExpansion
{\textstyle\int_{t}^{T}}
%EndExpansion
e^{\beta s}(\left\vert \bar{Y}(s)\right\vert ^{2}+\left\vert \bar
{Z}(s)\right\vert ^{2}+%
%TCIMACRO{\tint _{E}}%
%BeginExpansion
{\textstyle\int_{E}}
%EndExpansion
\left\vert \bar{K}(s,\zeta)\right\vert ^{2}\nu(d\zeta))ds|\mathcal{F}_{t}]\\
\leq\tfrac{1}{\rho}\mathbb{E}[%
%TCIMACRO{\tint _{t}^{T}}%
%BeginExpansion
{\textstyle\int_{t}^{T}}
%EndExpansion
e^{\beta s}%
%TCIMACRO{\tint _{-\delta}^{0}}%
%BeginExpansion
{\textstyle\int_{-\delta}^{0}}
%EndExpansion
(\left\vert \bar{U}(s+r)\right\vert ^{2}+\left\vert \bar{V}(s+r)\right\vert
^{2}+%
%TCIMACRO{\tint _{E}}%
%BeginExpansion
{\textstyle\int_{E}}
%EndExpansion
\left\vert \bar{Q}(s+r,\zeta)\right\vert ^{2}\nu(d\zeta))\mu(dr)ds|\mathcal{F}%
_{t}].
\end{array}
\label{nacira}%
\end{equation}
By changing the variables\textbf{\ }$v=s+r,$\textbf{\ }we have%
\begin{equation}%
\begin{array}
[c]{l}%
%TCIMACRO{\tint _{0}^{T}}%
%BeginExpansion
{\textstyle\int_{0}^{T}}
%EndExpansion
e^{\beta s}%
%TCIMACRO{\tint _{-\delta}^{0}}%
%BeginExpansion
{\textstyle\int_{-\delta}^{0}}
%EndExpansion
(\left\vert \bar{U}(s+r)\right\vert ^{2}+\left\vert \bar{V}(s+r)\right\vert
^{2}+%
%TCIMACRO{\tint _{E}}%
%BeginExpansion
{\textstyle\int_{E}}
%EndExpansion
\left\vert \bar{Q}(s+r,\zeta)\right\vert ^{2}\nu(d\zeta))\mu(dr)ds\\
=%
%TCIMACRO{\tint _{-\delta}^{0}}%
%BeginExpansion
{\textstyle\int_{-\delta}^{0}}
%EndExpansion
e^{-\beta r}%
%TCIMACRO{\tint _{0}^{T}}%
%BeginExpansion
{\textstyle\int_{0}^{T}}
%EndExpansion
e^{\beta(s+r)}(\left\vert \bar{U}(s+r)\right\vert ^{2}+\left\vert \bar
{V}(s+r)\right\vert ^{2}+%
%TCIMACRO{\tint _{E}}%
%BeginExpansion
{\textstyle\int_{E}}
%EndExpansion
\left\vert \bar{Q}(s+r,\zeta)\right\vert ^{2}\nu(d\zeta))ds\mu(dr)\\
=%
%TCIMACRO{\tint _{-\delta}^{0}}%
%BeginExpansion
{\textstyle\int_{-\delta}^{0}}
%EndExpansion
e^{-\beta r}%
%TCIMACRO{\tint _{r}^{T+r}}%
%BeginExpansion
{\textstyle\int_{r}^{T+r}}
%EndExpansion
e^{\beta v}(\left\vert \bar{U}(v)\right\vert ^{2}+\left\vert \bar
{V}(v)\right\vert ^{2}+%
%TCIMACRO{\tint _{E}}%
%BeginExpansion
{\textstyle\int_{E}}
%EndExpansion
\left\vert \bar{Q}(v,\zeta)\right\vert ^{2}\nu(d\zeta))dv\mu(dr)\\
\leq%
%TCIMACRO{\tint _{-\delta}^{0}}%
%BeginExpansion
{\textstyle\int_{-\delta}^{0}}
%EndExpansion
e^{-\beta r}\mu(dr)(\left\vert r\right\vert \left\vert \bar{U}(0)\right\vert
^{2}+%
%TCIMACRO{\tint _{0}^{T}}%
%BeginExpansion
{\textstyle\int_{0}^{T}}
%EndExpansion
e^{\beta v}(\left\vert \bar{U}(v)\right\vert ^{2}+\left\vert \bar
{V}(v)\right\vert ^{2}+%
%TCIMACRO{\tint _{E}}%
%BeginExpansion
{\textstyle\int_{E}}
%EndExpansion
\left\vert \bar{Q}(v,\zeta)\right\vert ^{2}\nu(d\zeta))dv)\mathbf{.}%
\end{array}
\label{elin}%
\end{equation}
Combining\ $\left(  \ref{elin}\right)  $ with\ $\left(  \ref{nacira}\right)  $
and taking\ $t=0$ and taking conditional expectation,\ we obtain%
\[%
\begin{array}
[c]{l}%
\mathbb{E}[\left\vert \bar{Y}(0)\right\vert ^{2}]+\mathbb{E}[%
%TCIMACRO{\tint _{0}^{T}}%
%BeginExpansion
{\textstyle\int_{0}^{T}}
%EndExpansion
e^{\beta s}(\left\vert \bar{Y}(s)\right\vert ^{2}+\left\vert \bar
{Z}(s)\right\vert ^{2}+%
%TCIMACRO{\tint _{E}}%
%BeginExpansion
{\textstyle\int_{E}}
%EndExpansion
\left\vert \bar{K}(s,\zeta)\right\vert ^{2}\nu(d\zeta))ds]\\
\leq\tfrac{1}{\rho}%
%TCIMACRO{\tint _{-\delta}^{0}}%
%BeginExpansion
{\textstyle\int_{-\delta}^{0}}
%EndExpansion
e^{-\beta r}\mu(dr)\mathbb{E}[\left\vert \bar{U}(0)\right\vert ^{2}+%
%TCIMACRO{\tint _{0}^{T}}%
%BeginExpansion
{\textstyle\int_{0}^{T}}
%EndExpansion
e^{\beta s}(\left\vert \bar{U}(s)\right\vert ^{2}+\left\vert \bar
{V}(s)\right\vert ^{2}+%
%TCIMACRO{\tint _{E}}%
%BeginExpansion
{\textstyle\int_{E}}
%EndExpansion
\left\vert \bar{Q}(s,\zeta)\right\vert ^{2}\nu(d\zeta))ds].
\end{array}
\]
for\ $\delta\in\left(  0,1\right)  .$ As
\[
\tfrac{1}{\rho}%
%TCIMACRO{\tint _{-\delta}^{0}}%
%BeginExpansion
{\textstyle\int_{-\delta}^{0}}
%EndExpansion
e^{-\beta r}\mu(dr)\rightarrow\tfrac{1}{\rho}\mu(\{0\})<1,\text{ as }%
\delta\rightarrow0,
\]
there is some $\delta_{\rho}>0$ such that for all $\delta\in\left(
0,\delta_{\rho}\right)  :\tfrac{1}{\rho}%
%TCIMACRO{\tint _{-\delta}^{0}}%
%BeginExpansion
{\textstyle\int_{-\delta}^{0}}
%EndExpansion
e^{-\beta r}\mu(dr)<1,$ i.e., $\Phi:(\left(  L^{2}(\mathcal{F}_{0}%
)\times\mathbf{L}_{T}^{2}\right)  \times\mathbf{L}_{T}^{2}\times\mathbf{H}%
_{T}^{2},\left\Vert \cdot\right\Vert _{\beta})$\textbf{$\rightarrow$}$(\left(
L^{2}(\mathcal{F}_{0})\times\mathbf{L}_{T}^{2}\right)  \times\mathbf{L}%
_{T}^{2}\times\mathbf{H}_{T}^{2},\left\Vert \cdot\right\Vert _{\beta}%
)$\textbf{\ }has a unique fixed point:\textbf{\ }$\left(  \left(
Y(0),Y\right)  ,Z,K\right)  $ $\in$ $(\left(  L^{2}(\mathcal{F}_{0}%
)\times\mathbf{L}_{T}^{2}\right)  \times\mathbf{L}_{T}^{2}\times\mathbf{H}%
_{T}^{2},\left\Vert \cdot\right\Vert _{\beta})$\textbf{$.\qquad\qquad\qquad$%
}$\square$


\begin{thebibliography}{99}                                                                                               %


\bibitem {AO3}Agram, N., \& Øksendal, B. (2016). Model uncertainty stochastic
mean-field control. arXiv preprint arXiv:1611.01385.

\bibitem {AO2}Agram, N., \& Øksendal, B. (2017). Stochastic Control of Memory
Mean-Field Processes. Applied Mathematics \& Optimization, 1-24.

\bibitem {A}Agram, N. (2016). Stochastic optimal control of McKean-Vlasov
equations with anticipating law. arXiv preprint arXiv:1604.03582.

\bibitem {AR}Agram, N., \& Røse, E.E. (2017). Optimal control of
forward-backward mean-field stochastic delay systems. Afr. Mat. https://doi.org/10.1007/s13370-017-0532-6.

\bibitem {AR2}Agram, N., \& Røse, E.E. (2017). Mean-field delayed BSDEs in
finite and infinite horizon. arXiv preprint arXiv:1509.08777.

\bibitem {BBP}Barles, G., Buckdahn, R., \& Pardoux, E. (1997). Backward
stochastic differential equations and integral-partial differential equations.
Stochastics: An International Journal of Probability and Stochastic Processes,
60(1-2), 57-83.

\bibitem {b}Bismut, J. M. (1978). An introductory approach to duality in
optimal stochastic control. SIAM review, 20(1), 62-78.

\bibitem {B1}Buckdahn, R., Li, J., \& Peng, S. (2009). Mean-field backward
stochastic differential equations and related partial differential equations.
Stochastic Processes and their Applications, 119(10), 3133-3154.

\bibitem {CD}Carmona, R., \& Delarue, F. (2015). Forward--backward stochastic
differential equations and controlled McKean--Vlasov dynamics. The Annals of
Probability, 43(5), 2647-2700.

\bibitem {di}Delong, \L ., \& Imkeller, P. (2010). Backward stochastic
differential equations with time delayed generators---results and
counterexamples. The Annals of Applied Probability, 20(4), 1512-1536.

\bibitem {d}Delong, L. (2010). Applications of time-delayed backward
stochastic differential equations to pricing, hedging and portfolio
management. arXiv preprint arXiv:1005.4417.

\bibitem {dbo}Delong, \L . (2013). Backward stochastic differential equations
with jumps and their actuarial and financial applications. London: Springer.

\bibitem {DI2}Delong, \L ., \& Imkeller, P. (2010). On Malliavin's
differentiability of BSDEs with time delayed generators driven by Brownian
motions and Poisson random measures. Stochastic Processes and their
Applications, 120(9), 1748-1775.

\bibitem {EPQ}El Karoui, N., Peng, S., \& Quenez, M. C. (1997). Backward
stochastic differential equations in finance. Mathematical finance, 7(1), 1-71.

\bibitem {ml}Ma, H., \& Liu, B. (2017). Infinite horizon optimal control
problem of mean-field backward stochastic delay differential equation under
partial information. European Journal of Control.

\bibitem {SEM}Mohammed, S. E. A.: Stochastic differential equations with
memory: Theory, examples and applications. Stochastic analysis and related
topics VI. The Geilo Workshop, 1996, Progress in Probability, Birkhauser, 1998.

\bibitem {OS}Øksendal, B., \& Sulem, A. (2015). Risk minimization in financial
markets modeled by Itô-Lévy processes. Afrika Matematika, 26(5-6), 939-979.

\bibitem {SQ}Quenez, M. C., \& Sulem, A. (2013). BSDEs with jumps,
optimization and applications to dynamic risk measures. Stochastic Processes
and their Applications, 123(8), 3328-3357.

\bibitem {PP}Pardoux, E., \& Peng, S. (1990). Adapted solution of a backward
stochastic differential equation. Systems \& Control Letters, 14(1), 55-61.

\bibitem {P}Pardoux, É. (1999). BSDEs, weak convergence and homogenization of
semilinear PDEs. NATO ASI Series C Mathematical and Physical Sciences-Advanced
Study Institute, 528, 503-550.

\bibitem {R}Royer, M. (2006). Backward stochastic differential equations with
jumps and related non-linear expectations. Stochastic processes and their
applications, 116(10), 1358-1376.

\bibitem {TL}Tang, S., \& Li, X. (1994). Necessary conditions for optimal
control of stochastic systems with random jumps. SIAM Journal on Control and
Optimization, 32(5), 1447-1475.
\end{thebibliography}
\end{document}